\documentclass{elsart}
\usepackage{amsmath,amscd, amsfonts}
\def\bg{\operatorname{beg}}
\def\ed{\operatorname{end}}
\def\gen{\operatorname{gen}}
\def\hdeg{\operatorname{hdeg}}
\def\deg{\operatorname{deg}}
\def\adeg{\operatorname{adeg}}
\def\reg{\operatorname{reg}}
\def\Tor{\operatorname{Tor}}
\def\Ext{\operatorname{Ext}}
\def\Ass{\operatorname{Ass}}
\def\gin{\operatorname{gin}}
\def\Gin{\operatorname{Gin}}
\def\Hom{\operatorname{Hom}}
\def\dept{\operatorname{depth}}
\def\sk{\smallskip}

\def\Zset{{\mathbb Z}}
\def\mm{{\mathfrak m}}
\def\pfrak{{\mathfrak p}}
\def\qfrak{{\mathfrak q}}
\def\Hcal{{\mathcal H}}
\def\Mcal{{\mathcal M}}
\def\Mbar{\overline{M}}

\begin{document}
\begin{frontmatter}

\title{ Castelnuovo-Mumford regularity\\ of canonical and deficiency modules}
\thanks{Both authors were supported by the Academy of Finland, Project 48556.
The first author was also supported  in part by the National Basic Research
Program (Vietnam).}
\author{ L\^e Tu\^an Hoa}
\address{ Institute of Mathematics Hanoi, 18 Hoang Quoc Viet Road, 10307 Hanoi, 
Vietnam}
\ead{lthoa@math.ac.vn}
\author{Eero Hyry}
\address{ 
Department of Mathematics, University of Helsinki, 
PL 68 (Gustaf H\"allstr\"omin katu 2b)
FIN--00014 HELSINGIN YLIOPISTO,
Finland }
\ead{eero.hyry@helsinki.fi}
\begin{abstract}  We give two kinds of bounds for the Castelnuovo-Mumford
regularity of the canonical module and the deficiency modules of a ring, 
respectively in terms 
of the homological degree and the Castelnuovo-Mumford regularity of the original 
ring. 
\end{abstract}
 \maketitle

\end{frontmatter}

\section*{Introduction}\smallskip

Let $R= k[x_1,...,x_n]$ be a standard graded polynomial ring, and let $M$ be a
finitely generated graded $R$-module of dimension $d$. The canonical module
$K^d(M) = \Ext^{n-d}_R(M,R)(-n)$ of $M$ - originally introduced by
Grothendieck - plays an important role in  Commutative Algebra and Algebraic
Geometry (see, e.g., \cite{BH}). It is natural to ask whether one can bound
the Castelnuovo-Mumford regularity $\reg(K^d(M))$ of $K^d(M)$ in terms of
other invariants of $M$. Besides the canonical module we are also
interested in a similar problem for all the deficiency modules $K^i(M) =
\Ext^{n-i}_R(M,R)(-n),\ i<d$. These modules were defined in \cite{Sc1},
Section 3.1, and can be considered as a measure for the deviation of $M$ being
a Cohen-Macaulay module.  Moreover, even in a rather simple case there is a
close relationship between $\reg(K^d(M))$ and $\reg(K^i(M)),\ i<d$ (see
\cite{Sc1}, Corollary 3.1.3). One could say that the
Castelnuovo-Mumford regularity $\reg(M)$ controls positive components of all
local cohomology modules $H^i_{\mm}(M)$ of $M$: they vanish above the level
$\reg(M)$. Although the negative components $H^i_{\mm}(M)_j$ do not
necessarily vanish, the function $\ell(H^i_{\mm}(M)_j)$ becomes a polynomial
for $j< -\reg(K^i(M))$. In this sense $\reg(K^i(M))$ controls the behavior of
$\ell(H^i_{\mm}(M)_j)$ in negative components. Note that in a series of papers
M. Brodmann and others have considered the problem when $\ell(H^i_{\mm}(M)_j)$
becomes a polynomial
 (see, e.g., \cite{BL}, \cite{BMM}). In fact, a result of \cite{BMM} will play
 an important role in our investigation.

 In this paper we will give two kinds of bounds for $\reg(K^i(M))$. In Section
 \ref{B} we will show that one can use the homological degree to bound
 $\reg(K^i(M))$. The homological degree was introduced by W. Vasconcelos
 \cite{Va1}, and one can use it to bound $\reg(M)$ (see \cite{DGV}, Theorem 2.4
 and \cite{Na}, Theorem 3.1). For the ring case the bound has a simple form: $\reg(S)
 < \hdeg(S)$, where $S$ is a quotient ring of $R$. Our Theorem
 \ref{B4} says that $\reg(K^i(S)) \leq d\cdot \hdeg(S)$ for all $i$. Thus this result
 complements the relationship between the Castelnuovo-Mumford regularity and
 the homological degree.

 The core of the paper is Section \ref{C}. Here  we restrict ourselves to the
 case of rings. We will then prove that one can bound $\reg(K^i(S))$ in terms of
 $\reg(S)$ (see Theorem \ref{C1}). Although the bounds are huge numbers, they
 show that the Castelnuovo-Mumford regularity also controls the behavior of
 all  local cohomology modules in negative components (in the above mentioned
 sense). This is a new meaning for the  Castelnuovo-Mumford regularity. For
example, our study gives the following consequence

\vskip0.3cm
 \noindent {\bf Corollary \ref{C5}}. {\it Denote by $\Hcal_{n,i,r}$ the set of  numerical
functions $h:\ \Zset \rightarrow \Zset$ such that there exists a homogeneous
ideal $I \subset R= k[x_1,...,x_n]$ satisfying the following conditions
\begin{itemize}
\item[(i)] $\reg I \leq r$,
\item[(ii)] $\ell(H^i_{\mm}(R/I)_t) = h(t)$ for all $t\in \Zset$.
\end{itemize}
Then for fixed numbers $n,i,r$ the set $\Hcal_{n,i,r}$ has only finitely many
elements.}
 \vskip0.3cm

In the last section \ref{D} we will examine some cases where $\reg(K^i(M))$ can be
bounded by a small number. As one can expect, each of them is of a very special
type. Section \ref{A}, where we collect some
results on the Castelnuovo-Mumford regularity, is of preparatory character.

\section{Preliminaries}\smallskip \label{A}

In this section we recall some basis facts on the Castelnuovo-Mumford
regularity. Throughout the paper let $R= k[x_1,...,x_n]$ be a standard graded
polynomial ring, where $k$ is an infinite field, and let $\mm =
(x_1,...,x_n)$. For an arbitrary graded $R$-module $N$, put
$$\bg(N) = \inf\{i \in \Zset|\ [N]_i \neq 0\}, $$
and
$$\ed(N) = \sup\{i \in \Zset|\ [N]_i \neq 0\}.$$

\noindent (We assume $\bg(N) = + \infty$ and $\ed(N) = -\infty$ if $N=0$.)

\begin{defn}\label{A1}{\rm Let $M$ be a finitely generated $R$-module.
The number
$$\reg (M) =  \max\{i + \ed(H^i_{\mm}(M))|\ i\geq  0\} $$
 is called the Castelnuovo-Mumford regularity of $M$.}
\end{defn}

Note that if $I\subset R$ is a nonzero homogeneous ideal, then
$$\reg (I) = \reg (R/I) +1.$$
We also consider the number
$$\reg_1(M) =  \max\{i + \ed(H^i_{\mm}(M))|\ i\geq 1\},$$
which is sometimes called the Castelnuovo-Mumford regularity at level one. The
definition immediately gives
\begin{equation} \label{EA1}
\reg (M) = \max\{\reg_1(M),\ \ed(H^0_{\mm}(M))\}.
\end{equation}
The following result is the starting point for the investigation of the
Castelnuovo-Mumford regularity.

\begin{lem}\label{A2}{\rm (\cite{EG}, Proposition 1.1 and Theorem 1.2)}
$$\reg (M )=\max\{\ed(\Tor^R_i(k,M)) -i|\ i\geq 0\}.$$
\end{lem}

The long exact sequence of local cohomology arising from a short exact
sequence of modules gives:

\begin{lem}\label{A3}{\rm (\cite{E}, Corollary 20.19)} Let
$$0 \rightarrow A \rightarrow B \rightarrow C \rightarrow 0$$
be an exact sequence of graded $R$-modules. Then
\begin{itemize}
\item[(i)] $\reg (B) \leq \max\{\reg (A),\ \reg (C)\}$,
\item[(ii)] $\reg (A) \leq \max\{\reg (B),\ \reg (C) +1\}$.
\end{itemize}
\end{lem}

Recall that a homogeneous element $x\in \mm$ is called an $M$-filter regular
if
$$x \not\in \pfrak\ \ \text{for \ all}\ \ \pfrak \in (\Ass M)\setminus \{
\mm\}.$$ This is equivalent to the condition that the module $0:_M x$ is  of
finite length. Since $k$ is assumed to be infinite, there always exists a
filter regular element with respect to a finite number of finitely generated
modules.

Let $x$ be a linear $M$-filter regular element. Then
$$H^i_{\mm}(M/0:_M x) \cong H^i_{\mm}(M) \ \ \text{for \ all}\ \ i\geq 1.$$
Hence the short exact sequence induced by multiplication by $x$
$$ 0\rightarrow (M/0:_M x)(-1) \overset{\cdot x}\rightarrow M \rightarrow M/xM
\rightarrow 0$$ provides the exact sequence
$$\begin{array}{ll}
0 & \rightarrow (0:_M x)_{j-1} \rightarrow H^0_{\mm}(M)_{j-1} \rightarrow
H^0_{\mm}(M)_j \rightarrow H^0_{\mm}(M/xM)_j \rightarrow \cdots \\
& \cdots \rightarrow H^i_{\mm}(M)_j \rightarrow H^i_{\mm}(M/xM)_j \rightarrow
H^{i+1}_{\mm}(M)_{j-1} \rightarrow H^{i+1}_{\mm}(M)_j \rightarrow \cdots
\end{array}$$
From this one can get (see \cite{E}, Proposition 20.20 and \cite{HH}, Lemma
2):

\begin{lem}\label{A5} Let $x$ be a linear $M$-filter
regular element. Then
$$ \reg_1(M) \leq \reg (M/x M) \leq \reg M.$$
\end{lem}

Finally, let us recall the notion of the regularity index (of a Hilbert
function). In the literature it also appears under different names like the
$a$-invariant (see \cite{BH}, Definition 4.3.6 and Theorem 4.3.5,  and
\cite{Va2}, Section B.4) or the postulation number \cite{BMM}.

\begin{defn}\label{A6} {\rm Let $H_M(t)$ and $P_M(t)$ denote the Hilbert
function and the Hilbert polynomial of $M$, respectively. The number $$ri(M) =
\max\{j \in \Zset|\ H_M(j) \neq P_M(j)\}$$ is called the regularity index of $M$.}
\end{defn}

\begin{lem}\label{A7} Let $x$ be a linear $M$-filter
regular element. Then
\begin{itemize}
\item[(i)] {\rm (cf. \cite{E}, Proposition 20.20)} $\reg (M )= \max\{\reg (M/xM),\ 
\ed(H^0_{\mm}(M))\}$,
\item[(ii)] $\reg (M) = \max\{\reg (M/xM),\ ri(M)\}$,
\item[(iii)] If $M$ is a Cohen-Macaulay module of dimension $d$, then $\reg (M) = 
ri(M) + d$.
\end{itemize}
\end{lem}

\begin{pf} (i) This follows from Lemma \ref{A5} and (\ref{EA1}).

(ii) From the Grothendieck-Serre formula
\begin{equation}\label{EA3}
H_M(j) - P_M(j) = \sum_{i=0}^d (-1)^i \ell (H^i_{\mm}(M)_j),
\end{equation}
it follows that $\reg M \geq ri(M)$. By Lemma \ref{A5} we get
$$\reg (M ) \geq \max\{\reg (M/xM),\ ri(M)\}.$$
Let $j \geq \reg (M/xM)$. Since $\reg_1(M) \leq j$,  this yields by (\ref{EA3})
$$H_M(j) - P_M(j) = \ell (H^0_{\mm}(M)_j).$$
Hence
$$\ed (H^0_{\mm}(M)) \leq \max\{\reg (M/xM),\ ri(M)\}.$$
Together with (i) we get
$$\reg (M) \leq \max\{\reg (M/xM),\ ri(M)\}.$$

(iii) This follows from (\ref{EA3}) and the fact that $H^i_{\mm}(M) = 0$ for
all $i< d$.
\end{pf}

\section{Comparison with homological degree}\smallskip \label{B}

From now on let $M$ be a $d$-dimensional finitely generated graded $R$-module.
The homological degree of a graded $R$-module $M$ was introduced by
Vasconcelos. It is defined recursively on the dimension as follows:

\begin{defn}\label{B1} {\rm (\cite{Va1} and \cite{Va2}, Definition 9.4.1)
The homological degree of $M$ is the number $$\hdeg (M) = \deg (M) +
\sum_{i=0}^{d-1}{d-1 \choose i} \hdeg(\Ext_R^{n+i+1-d}(M,R)).$$}
\end{defn}

Note that

(a) $\hdeg (M) \geq  \deg(M) $, and the equality holds if and only if $M$ is a
Cohen-Macaulay module.

(b) $\hdeg (M) = \hdeg(M/H^0_{\mm}(M)) + \ell(H^0_{\mm}(M))$.

Let $\gen (M)$ denote the maximal degree of elements in a minimal set of
homogeneous generators of $M$. That is,
$$\gen (M) = \ed(M/ \mm M).$$
It turns out that the homological degree gives an upper bound for the
Castelnuovo-Mumford regularity

\begin{lem}\label{B2} {\rm (\cite{DGV}, Theorem 2.4 and \cite{Na}, Theorem
3.1)}
$$\reg (M )\leq \gen (M) + \hdeg (M) -1.$$
\end{lem}

 Let
$$K^i(M) = \Ext^{n-i}_R(M,R)(-n).$$
The module $K^d(M)$ is the canonical module of $M$. Following Schenzel
(\cite{Sc1}, Section 3.1) we call the modules $K^i(M),\ i< d,$ as the deficiency modules 
of $M$.
Note that $K^i(M) = 0$ for $i<0$ and $i>d$. All the modules $K^i(M)$ are finitely
generated, and by \cite{Sc1}, Section 3.1 (see Lemma 3.1.1 and page 63) we
have:
$$\begin{array}{l}
\dim K^i(M) \leq i \ \text{for}\ i<d,\\
  \dim K^d(M) = d, \ \text{and}\\
 \dept (K^d(M)) \geq \min\{ 2, \dim M\}. 
\end{array}$$
By the local duality theorem (see, e.g., \cite{BH}, Theorem 3.6.19), there are
the following canonical isomorphisms of graded modules
\begin{equation}\label{EB2}
K^i(M) \cong \Hom_k (H^{i}_{\mm}(M), k).
\end{equation}
 From this and Lemma \ref{A7} (ii) we obtain that
$$\ell (H^{i}_{\mm}(M)_t)= P_{K^i(M))}(-t)\ \ \text{for\ all} \ t< -
\reg(K^i(M)).$$

Inspired by Lemma \ref{B2} it is natural to ask whether  one can use the
homological degree to bound the Castelnuovo-Mumford regularity of $K^i(M)$,
too? The following theorem, which is the main result of this section, answers this question affirmatively.

\begin{thm}\label{B4} For all $i\leq d$ we have
$$\reg (K^i(M)) \leq d[\hdeg (M) - \deg (M)] - \bg(M)  + i.$$
\end{thm}

Note that when $M$ is a Cohen-Macaulay module, $K^d(M)$ is also a
Cohen-Macaulay module. It was shown in \cite{HaH}, Proposition 2.3 that
\begin{equation}\label{EB2b}
\reg (K^d(M)) = d - \bg(M).
\end{equation}
This easily follows from Lemma \ref{A7} (iii) and the Grothendieck-Serre
formula (\ref{EA3}) applied to $K^d(M)$, or from the duality. Thus in this
case we have the equality in (ii) of the above theorem.

In order to prove Theorem \ref{B4} we need some auxiliary results.

\begin{lem}\label{B7} {\rm (\cite{Sc3}, Proposition 2.4)} Let
$x$ be a linear $M$-filter regular element. Then there are short exact
sequences of graded modules
$$ 0 \rightarrow (K^{i+1}(M) / x K^{i+1}(M))(1) \rightarrow K^i(M/xM)
\rightarrow 0:_{K^i(M)} x \rightarrow 0,$$ for all integers $i\geq 0$.
\end{lem}

For short, in the proof we often use the following notation
$$K^i := K^i(M).$$
\vskip0.3cm

\begin{lem}\label{B7bn} $\reg(K^0(M)) \leq - \bg(M)$. \end{lem}

\begin{pf} Note that $H^0_{\mm}(M) \subseteq M$ is a submodule of finite length.
Hence, by (\ref{EB2}), we have
$$ \reg(K^0) = - \bg( H^0_{\mm}(M)) \leq - \bg(M).$$
\hfill $\square$ \end{pf}

In the sequel we always assume that $x$ is a generic linear element by which we
mean that $x$ is filter regular with respect to $M$, all the modules $K^i(M)$ and
all the iterated deficiency modules in the sense of \cite{Va1}, Definition 2.12.
Since this is a finite collection of modules, such an element always exists.

 \vskip0.3cm

\begin{lem}\label{B7cn} Assume $\dept(M) > 0$ and $1 \leq i < d$. Then
$$\reg (K^i(M)) \leq -\bg(M) + \sum_{j=1}^i {d  \choose j} \hdeg (K^j(M))  +
i.$$
\end{lem}

\begin{pf} Let $x\in R$ be a generic linear element and $j\geq 0$. By Lemma \ref{B7} 
there
is an exact sequence
$$ 0 \rightarrow (K^{j+1} / x K^{j+1})(1) \rightarrow K^j(M/xM)
\rightarrow 0:_{K^j} x \rightarrow 0.$$
Taking the tensor product with $k$ we
get the exact sequence
$$K^j(M/xM)/\mm K^j(M/xM) \leftarrow (K^{j+1} / \mm K^{j+1})(1)
\leftarrow \Tor^R_1(k, 0:_{K^j} x).$$
This implies that
$$\begin{array}{ll}
\gen(K^{j+1}) & =  \ed(K^{j+1} / \mm K^{j+1}) \\
& \leq  \max\{ \gen(K^j(M/xM)),\ \ed(\Tor^R_1(k, 0:_{K^j}
x)) \} + 1.
\end{array}$$
Since $0:_{K^j} x$ is of finite length,
$$0:_{K^j} x \subseteq H^0_{\mm}(K^j).$$
Hence, by Lemma \ref{A2},
$$\ed(\Tor^R_1(k, 0:_{K^j} x)) -1 \leq \reg (0:_{K^j} x ) \leq \ed
(H^0_{\mm}(K^j)) \le \reg(K^j).$$
Combining this with the fact that 
$$\gen(K^j(M/xM)) \leq \reg (K^j(M/xM))$$
(look again at Lemma \ref{A2}),
we get
\begin{eqnarray} \label{EB5}
\gen(K^{j+1}) & \leq & \max\{ \reg(K^j(M/xM)),\ \reg (K^j)  + 1\} + 1 \nonumber \\
 & \leq& \max\{ \reg(K^j(M/xM)) +  1,\ \reg (K^j)  +2\} .
\end{eqnarray}
Note that
$$\bg(M/xM) \geq \bg(M).$$
We now prove the claim by induction on $i$. Let $i=1$.  An application of  (\ref{EB5}) to 
the case $j=0$ 
together with Lemma \ref{B7bn} yields
$$\gen(K^1) \leq \max\{ - \bg(M/xM),\ -\bg(M) \}  +2 = -\bg(M) +2.$$
By Lemma \ref{B2}, we then get
$$\begin{array}{ll}
 \reg(K^1) & \leq \gen(K^1) + \hdeg(K^1) -1 \leq \hdeg(K^1)-\bg(M) + 1\\
& \leq d\cdot \hdeg(K^1)-\bg(M)  + 1.
\end{array}$$
Thus the claim holds for $K^1$.  

Let $2\leq i\leq d-1$. By the induction hypothesis we have
\begin{equation}\label{EB5aa}
\reg (K^{i-1} ) \leq -\bg(M) + \sum_{j=1}^{i-1} {d \choose j} \hdeg (K^j) +
i - 1.
\end{equation}
For a Noetherian graded module $N$ over $S$, let $\overline{N}$ denote the module $N/H^0_\mm(N)$. Note that $\dept(\overline{N}) >0$ if $\dim N >0$, and for all $j>0$ we have
\begin{equation}\label{EB5ab}
K^j(N) \cong K^j(\overline{N}).
\end{equation}
Since  $\dim \overline{M/xM} = d-1$ and 
 $0< i-1 < d-1$, again by the  induction hypothesis applied to $\overline{M/xM} $, the following holds
$$\begin{array}{l}
\reg (K^{i-1} (M/xM)) =  \reg (K^{i-1} (\overline{M/xM} ))  \\
\quad \leq  -\bg(\overline{M/xM} ) + \displaystyle{ \sum_{j=1}^{i-1} {d-1 \choose j} \hdeg 
(K^j(\overline{M/xM} )) } + i -1 \\
\quad \leq   -\bg(M/xM) + \displaystyle{  \sum_{j=1}^{i-1} {d-1 \choose j} \hdeg 
(K^j(M/xM)) } + i - 1.
\end{array}$$
Since $\dept (M) >0$,  we have by the inequality (10) in \cite{Va1}
$$\hdeg (K^j(M/ xM)) \leq \hdeg K^j + \hdeg K^{j+1}.$$
So
$$\begin{array}{l}
\reg (K^{i-1} (M/xM))  \leq\\
\quad \leq -\bg(M) + \displaystyle{\sum_{j=1}^{i-1}} {d-1 \choose 
j} (\hdeg (K^j) +\hdeg(K^{j+1})) + i - 1\\
\quad \leq  -\bg(M) + \displaystyle{\sum_{j=1}^{i-1}} {d \choose j} \hdeg (K^j) + {d-1
\choose i-1} \hdeg (K^{i}) + i - 1. \hskip0.5cm {\rm (*)}
\end{array}$$
By (\ref{EB5}) and (\ref{EB5aa}) this  yields
$$\gen(K^i) \leq  -\bg(M) + \displaystyle{\sum_{j=1}^{i-1}} {d \choose j} \hdeg (K^j) + {d-1
\choose i-1} \hdeg (K^{i}) + i +1. \hskip0.5cm {\rm (**)}$$
Hence, by Lemma \ref{B2}, we then get 
$$\begin{array}{ll}
\reg (K^{i}) &  \leq \gen(K^{i}) + \hdeg(K^{i}) - 1\\
 &  \leq -\bg(M) + \displaystyle{ \sum_{j=1}^{i-1}}{d \choose j} \hdeg (K^j)
+ {d-1 \choose i-1 } \hdeg (K^{i}) + i +1\\
& \quad + \hdeg (K^{i}) - 1\\
&  \leq -\bg(M) + \displaystyle{ \sum_{j=1}^{i-1}}{d \choose j} \hdeg (K^j)
+ {d-1 \choose i-1 } \hdeg (K^{i}) + i \\
& \quad +  \displaystyle{ {d-1 \choose i }}  \hdeg (K^{i}) \\
& =  \displaystyle{\sum_{j=1}^{i}} {d \choose j} \hdeg (K^j)- \bg(M) + i.
\end{array}$$
Lemma \ref{B7cn}  is thus completely proved. \hfill $\square$ 
\end{pf}

\noindent {\it Remark. } Let us take an extra look at  the case $i=d-1$, where $d\ge 2$. If the equality holds, i.e. 
$$\reg (K^{d-1}(M)) = -\bg(M) + \sum_{j=1}^{d-1} {d  \choose j} \hdeg (K^j(M))  +
d-1,$$
then we must have the equality in (**), too. Using  (*), (\ref{EB5aa}) and (\ref{EB5}), this yields $\hdeg(K^{d-1}(M))= 0$, or equivalently $ K^{d-1}(M) =0$.
\vskip0.3cm

\noindent {\bf PROOF OF THEOREM \ref{B4}.}  Since $\hdeg(M) \geq \deg(M)$, by 
Lemma \ref{B7bn} we may assume that $d\geq 1$ and $i\geq 1$. Let $\Mbar = M/H^0_\mm(M)$. 

(i)  First consider the case $1\leq i <d$. Since $\dept (\Mbar) >0$, the formula of $\hdeg(\Mbar)$ in 
Definition \ref{B1}  can be rewritten as follows 
$$\begin{array}{ll}
\hdeg(\Mbar) & = \deg(\Mbar) +  \displaystyle{  \sum_{j=1}^{d-1} {d-1 \choose j} } \hdeg 
(K^{d-j-1}(\Mbar))\\
& = \deg(\Mbar) +  \displaystyle{ \sum_{j=1}^{d-1} {d-1 \choose j} } \hdeg (K^{j}(\Mbar))\\
& \geq \deg(\Mbar)+  \displaystyle{  \sum_{j=1}^i {d-1 \choose j}} \hdeg (K^j(\Mbar)) \\
&\geq \deg(\Mbar)+  \displaystyle{ \frac{1}{d} \sum_{j=1}^i {d \choose j}} \hdeg (K^j(\Mbar)) .
\end{array}$$
Consequently, Lemma \ref{B7cn} and (\ref{EB5ab}) give 
$$\reg (K^i) = \reg(K^i(\Mbar)) \leq d(\hdeg (\Mbar) - \deg (\Mbar) )- \bg(\Mbar) + i.$$
Since $\hdeg(\Mbar) \leq \hdeg(M),\ \deg(\Mbar) = \deg(M)$ and $\bg(\Mbar) \geq 
\bg(M)$, the above inequality  yields
$$\reg (K^i) \leq d(\hdeg (M) - \deg M) - \bg(M)  + i.$$
Thus (i) is proved. \vskip0.3cm

(ii)  We now prove for the claim for $\reg (K^d(M)$. We  do induction on $d$.

If $d=1$, then $\Mbar$ is a Cohen-Macaulay module. By (\ref{EB2b}) we have
$$\begin{array}{ll}
\reg(K^1) & = \reg(K^1(\Mbar)) = 1 - \bg(\Mbar) \leq 1-\bg(M) \\
& \leq  1 + \hdeg (M) - \deg M - \bg(M) .
\end{array}$$
Let $d\geq 2$.  Let $x$ be a generic linear 
element.  Since $\dept K^d >0$,  one has by Lemma \ref{A7} (i) 
$$\reg(K^d) = \reg(K^d/ xK^d).$$
Using the short exact sequence
$$ 0 \rightarrow (K^d / x K^d)(1) \rightarrow K^{d-1}(M/xM)
\rightarrow 0:_{K^{d-1}} x \rightarrow 0,$$ and Lemma \ref{A3} (ii) we then
get
\begin{equation}\label{EB5b}
\reg(K^d) \leq \max\{ \reg(K^{d-1}(M/xM)),\ \reg (K^{d-1}) + 1\} + 1.
\end{equation}
If $K^{d-1} \neq 0$, then by Part (i)  and the remark after Lemma  \ref{B7cn} it already holds that
$ \reg K^{d-1} + 2  <   d(\hdeg(M) - \deg(M)) -\bg(M) + d +1$.
Hence 
$$ \reg K^{d-1} + 2  \le   d(\hdeg(M) - \deg(M)) -\bg(M) + d .$$
This inequality trivially holds if $K^{d-1}=0$. On the other hand, by the induction hypothesis
$$\begin{array}{ll}
\reg(K^{d-1}(M/xM))  & \leq (d-1)( \hdeg(M/xM) - \deg(M/xM)) -\bg(M/xM) + d-1\\
& \leq d( \hdeg(M/xM) - \deg(M)) -\bg(M) + d-1.
\end{array}$$
We now distinguish two cases:

\indent  $\bullet$  Assume $\dept M> 0$. By \cite{Va1}, Theorem 2.13,
we have
$$\hdeg (M/xM) \leq \hdeg(M).$$
Hence
$$\reg(K^{d-1}(M/xM)) +1 \leq d ( \hdeg(M) - \deg(M)) -\bg(M) + d.$$
 Summing up we obtain
$$\reg(K^d) \leq d( \hdeg(M) - \deg(M) )-\bg(M) +d .$$

\indent  $\bullet$  We now consider the case $\dept(M) =0$.  Since $\hdeg(\Mbar) \leq \hdeg(M),\ 
\deg(\Mbar) = \deg(M)$ and $\bg(\Mbar) \geq \bg(M)$,  by (\ref{EB5ab}) we get
$$\begin{array}{ll}
\reg(K^d) & = \reg(K^d(\Mbar)) \leq d( \hdeg(\Mbar) - \deg(\Mbar) )-\bg(\Mbar) +d \\
& \leq d( \hdeg(M) - \deg(M)) -\bg(M)  +d.
\end{array}$$
The proof of Theorem \ref{B4} is completed. \hfill $\square$ \vskip0.3cm

\begin{rem}{\rm Inspired by the homological degree Vasconcelos also
introduced a  class  of functions, called extended degree Deg($M$) (see
\cite{Va1} and \cite{Va2}, p. 263). This class contains $\hdeg(M)$. In fact,
Theorem 3.1 in \cite{Na} (see also \cite{DGV}, Theorem 2.4) states that 
$$\reg(M) \leq \gen(M) + \text{Deg}(M) -1.$$
It is interesting to ask whether one can replace $\hdeg(M)$ in Theorem
\ref{B4} by Deg($M$). Our method is not applicable in this case, because the
definition of an arbitrary extended degree does not explicitly contain the
information on $K^i(M)$.}
\end{rem}

\section{Castelnuovo-Mumford of a ring and its deficiency modules}\smallskip 
\label{C}

In this section we will consider a quotient ring $S=R/I$, and 
give a bound for $\reg(K^i(S)),\ i\leq d$, in terms of $\reg(S)$. We always
assume that $I$ is a {\it non-zero homogeneous ideal containing no linear form.}
 Note that  it is unclear whether one
can bound $\hdeg(S)$ in terms of $\reg(S)$. Therefore the following bound is
independent from that of Theorem \ref{B4}.

\begin{thm}\label{C1} Let $S=R/I$ be a quotient ring of a polynomial ring
$R = k[x_1,...,x_n]$ ($n\geq 2$)  modulo a homogeneous ideal $I \subseteq R$ as above. Then
$$\reg(K^i(S)) < \begin{cases} 4(\reg I)^{n-1} - 4 (\reg I)^{n-2} \ \
\text{if} \ \ i= 1, \\
(2 \reg I)^{n\cdots (n+i-1)2^{\frac{i(i-1)}{2}}} \ \ \text{if} \ \ i\geq 2.
\end{cases}$$
\end{thm}

In order to prove this theorem we need a result of M. Brodmann, C. Matteotti
and N. D. Minh \cite{BMM}. Following the notation there, we set
$$\begin{array}{ll}
h^i_S(t)&= \ell(H^i_{\mm}(S)_t) = H_{K^i(S)}(-t),\\
d^0_S(t) &= H_S(t) - h^0_S(t) + h^1_S(t),\\
d^i_S(t) &= h^{i+1}_S(t),\ i\geq 1.
\end{array}$$
Since $K^i(S)$ is a  finitely generated $R$-module, there is a polynomial
$q^i_S(t)$ such that
$$d^i_S(t) = q^i_S(t) \ \ \text{for} \ \ t\ll 0.$$
For $i\geq 0$, let
\begin{equation}\label{EC0}
\Delta_i = \sum_{j=0}^i {i \choose j}(d^j_S(-j) + | q^j_S(-j)|).
\end{equation}
Then Proposition 3.22 (c) of \cite{BMM} can be reformulated as follows

\begin{lem}\label{C2} For all $i\geq 1$ we have
$$ri(K^i(S)) \leq [2(1+\Delta_{i-1})]^{2^{i-1}} -2.$$
\end{lem}
\begin{pf} Set
$$\nu^i_S = \inf \{ t\in \Zset|\ d^i_S(t) \neq q^i_S(t) \}.$$
Proposition 3.22 (c) of \cite{BMM} states that
$$\nu^i_S \geq - [2(1+\Delta_{i})]^{2^i} + 2,$$
for all $i\geq 0$. Since $ri(K^i(S)) = - \nu^{i-1}_S$ for $i\geq 2$, the
assertion holds for $i\geq 2$. From the definition of $d^0_S(t)$ we also have
$$ri(K^1(S)) \leq \max \{0,\ - \nu^0_S\} \leq 2\Delta_0.$$
\hfill $\square$  \end{pf}

\begin{lem}\label{C3}
For $0 \leq i <d = \dim S$ and all $t\in \Zset$ we have
$$h^i_S(t)< (\reg I)^{n-i-1}{\reg(S) -t \choose i}.$$
\end{lem}

\begin{pf} This follows from \cite{Hoa}, Theorem 3.4 and Remark 3.5. \hfill $\square$  
\end{pf}

\noindent {\bf PROOF OF THEOREM \ref{C1}.} We divide the proof of Theorem \ref{C1} into proving
several claims. In the proof we simply write
$$K^i := K^i(S) \ \text{and}\ r := \reg(I) = \reg(S) +1.$$
\vskip0.3cm

 \noindent {\it CLAIM 1}. Let
$$\overline{K^{i+1}} = K^{i+1}/H^0_{\mm}(K^{i+1}),$$
and let $\alpha$ be an integer such that
$$\alpha \geq \max\{ \reg(K^i),\ \reg K^i(S/x S)\},$$
where $x$ is a generic linear element. Then
$$\reg(\overline{K^{i+1}}) \leq \reg(K^{i+1}/xK^{i+1}) \leq \alpha +2.$$

\begin{pf} Since $\reg(\overline{K^{i+1}})=\reg_1(K^{i+1})$, the first inequality 
holds by
 Lemma \ref{A5}. For the second inequality, by applying  Lemma \ref{A3}  to the exact sequence of
Lemma \ref{B7}
$$ 0 \rightarrow (K^{i+1} / x K^{i+1})(1) \rightarrow K^i(S/xS)
\rightarrow 0:_{K^i} x \rightarrow 0,$$
we get
$$\begin{array}{ll}
\reg(K^{i+1}/xK^{i+1}) &\leq 1 + \max\{ \reg(K^i(S/xS)),\ \reg(0:_{K^i} x)+1\}\\
&\leq 1 + \max\{ \reg(K^i(S/xS)),\ \ed (H^0_{\mm}(K^{i}))+1\}\\
&\leq 1 + \max\{ \reg(K^i(S/xS)),\ \reg (K^i)+1\}\\
&\leq \alpha + 2.
\end{array}$$
\hfill $\square$  
\end{pf}
\sk

\noindent {\it CLAIM 2}. $\reg(K^1) < 4(\reg I)^{n-1} - 4 (\reg I)^{n-2}$.

\begin{pf} Let $d=1$.  By (\ref{EB5ab}), $\reg(K^1) = \reg(K^1(\overline{S} ))$. Since 
$\overline{S}$  is a Cohen-Macaulay ring,  by  (\ref{EB2b})  we 
have $\reg(K^1) =1$. Since $I$ is a non-zero ideal and contains no linear form, $\reg I
\geq 2$. Hence the claim obviously holds in this case.

Let  $d\geq 2$. Since $H^0_{\mm}(S)$
is a submodule of finite length of $S$, $\reg (K^0) < 0$ for all $S$. Let $x$ be a generic 
linear element. By Claim
3 we get
\begin{equation}\label{EC2}
\reg(\overline{K^1}) \leq \reg(K^1/xK^1) \leq 2 + \max\{ \reg(K^0(S/xS)),\
\reg(K^0)\} \leq 1.
\end{equation}
Now we estimate $\Delta_0$. By Lemma \ref{C3},
$$\begin{array}{ll}
d^0_S(0) & = H_S(0) - h^0_S(0) + h^1_S(0) = 1 + h^1_S(0) \\
& < 1 + r^{n-2}\reg(S) = 1+ r^{n-1} - r^{n-2}.
\end{array}$$
Hence
$$d^0_S(0) \leq r^{n-1} - r^{n-2}.$$
  Note that
$$q_S^0(-t) = P_{K^1}(t) = P_{\overline{K^1}}(t).$$
Since $\reg (\overline{K^1}) \leq 1$,
$$h^1_{\overline{K^1}}(1) = 0.$$
As $\dim \overline{K^1} \leq 1$, applying the Grothendieck-Serre formula (see (\ref{EA3})) to $\overline{K^1}$
we get
$$P_{\overline{K^1}}(1) = H_{\overline{K^1}}(1).$$
Moreover, in this case $P_{\overline{K^1}}(t)$ is a constant. By Lemma
\ref{C3} this gives
$$\begin{array}{ll}
q^0_S(0) &= P_{\overline{K^1}}(0)= P_{\overline{K^1}}(1) = H_{\overline{K^1}}(1) 
\leq
H_{K^1}(1) \\
& = h^1_S(-1) < r^{n-2}\reg(S) =  r^{n-1} - r^{n-2}.
\end{array}$$
Putting all together we get
\begin{equation}\label{EC3}
\Delta_0 = d^0_S(0) + |q^0_S(0)| < 2r^{n-1} - 2r^{n-2}.
\end{equation}
From (\ref{EC2}) and Lemma \ref{C2} we can now conclude by  Lemma \ref{A7} (ii)  
that
$$\reg(K^1)\leq \max\{ 1, 2\Delta_0\}  < 4r^{n-1} - 4 r^{n-2}.$$
\hfill $\square$
\end{pf}

\sk

\noindent {\it CLAIM 3}. For $i\geq 1$ we have
$$\Delta_i < i\Delta_{i-1} + (\reg I)^{n-1} + |P_{\overline{K^{i+1}}}(-i)|.$$

\begin{pf} By (\ref{EC0}) we have
$$\begin{array}{ll}
\Delta_i & = \displaystyle{\sum_{j=0}^i} {i \choose j}(d^j_S(-j) + | q^j_S(-j)|) \\
 &=  \displaystyle{\sum_{j=0}^{i-1}} {i \choose j}(d^j_S(-j) + | q^j_S(-j)|) +
 h^{i+1}_S(-i) + |P_{K^{i+1}}(-i)|\\
 & \leq i\Delta_{i-1} + h^{i+1}_S(-i) + |P_{K^{i+1}}(-i)|.
 \end{array}$$
  By Lemma \ref{C3} we know that
$$h^{i+1}_S(-i) \leq r^{n-i-2}{r+i-1\choose i+1} < r^{n-i-2} r^{i+1} =
r^{n-1}.$$
Since $P_{K^{i+1}}(t) = P_{\overline{K^{i+1}}}(t)$, the claim
follows. \hfill $\square$
\end{pf}

\vskip0.3cm

\noindent {\it CLAIM 4}. Keep the notation and assumptions of Claim 1 with the additional assumption 
that $\alpha \geq 0$. 
For all $1\leq i < d-1$,
$$|P_{\overline{K^{i+1}}}(-i)| < \frac{1}{2} (\reg I)^{n-i-2}(\reg (I) +
\alpha + 2i + 1)^{2i+2}.$$

\begin{pf} By Claim  3, $\reg(\overline{K^{i+1}}) \leq \alpha+2$. Since
$\dept(\overline{K^{i+1}}) >0$ (if $\overline{K^{i+1}}\neq 0$), this implies
that
$$P_{\overline{K^{i+1}}}(t) = H_{\overline{K^{i+1}}}(t) \ \text{for\ all}\
t\geq \alpha +2.$$ In the case $P_{\overline{K^{i+1}}}(t) = 0$ there  is
nothing to prove. Assume that $\deg(P_{\overline{K^{i+1}}}(t)) = p \geq 0$.
Using the Lagrange's interpolation formula
$$\begin{array}{c}
 P_{\overline{K^{i+1}}}(t)  = \displaystyle{ \sum_{j=0}^p\frac{[t-(\alpha+2)] \cdots \widehat{[t -
(\alpha+2+j)]}\cdots [t - (\alpha+2 +p)]}{(j-0)(j-1)\cdots
\widehat{(j-j)}\cdots (j-p)} } \times \\
 P_{\overline{K^{i+1}}}(\alpha + 2+j),
\end{array}$$
 where $\hat{*}$ means that the corresponding term is omitted, we get
\begin{eqnarray} \label{EC3b}
P_{\overline{K^{i+1}}}(-i)& = \sum_{j=0}^p(-1)^j \frac{ (i+\alpha+2) \cdots
\widehat{(i+ \alpha+2+j)}\cdots  (i+\alpha+2 +p) }{|(j-0)(j-1)\cdots \widehat{(j-j)}\cdots (j-p)|} \times \nonumber \\
& \times H_{\overline{K^{i+1}}}(\alpha+ 2+j) . 
\end{eqnarray}
 Since $\dim  \overline{K^{i+1}} \leq i+1 $, $p\leq i+1$. By Lemma
\ref{C3} one has
$$\begin{array}{ll}
H_{ \overline{K^{i+1} }}(\alpha + 2+j) &\leq H_{K^{i+1}}(\alpha + 2+j) =
h^{i+1}_S( -(\alpha + 2+j)) \\
& \leq \displaystyle{ r^{n-i-2}{r-1 + \alpha + 2+j \choose i+1} } \\
& <  \displaystyle{ r^{n-i-2} \frac{(r+ \alpha + 1+ p)^{i+1}}{(i+1)!} }\\
&\leq \frac{1}{2} r^{n-i-2}(r+ \alpha + i+2)^{i+1}.
\end{array}$$
for all $j\leq p$.
Obviously
$$(i+\alpha+2) \cdots \widehat{(i+ \alpha+2+j)}\cdots  (i+\alpha+2 +p) \leq
(\alpha +i+2+p)^p \leq (\alpha+2i + 3)^{i+1}.$$
 Since $r \geq 2$,  the above estimations imply that all
numerators in (\ref{EC3b}) are strictly less than
$$ A:= \frac{1}{2} r^{n-i-2}(r+ \alpha + 2i+1)^{2i+2}.$$
All the denominators in the alternating sum (\ref{EC3b}) are bigger or equal
to $([\frac{p}{2}])^2$. There are at most $[\frac{p}{2}] +1$ terms with the same sign. 
This implies  that the sub-sum of all terms with the same sign in (\ref{EC3b}) has the absolute value less 
than $A$ if $p\geq 4$.  The same holds for $p\leq 3$ by a direct checking. Hence
$|P_{\overline{K^{i+1}}}(-i)| < A$. 

\hfill $\square$
\end{pf}

\sk

\noindent {\it CLAIM 5}. Assume that $d\geq 3$. Then
$$\Delta_1 <\frac{1}{2}(2 \reg(I))^{n(n+1)} - (\reg I)^n -n,$$
and
$$\reg(K^2) <(2 \reg(I))^{2n(n+1)}- 2(\reg I)^n -2n.$$

\begin{pf} By Claim 2
$$\reg(K^1) < 4 r^{n-1} - 4 r^{n-2} =: \alpha .$$
 Let $x$ be a generic linear element.  By Lemma \ref{A7} (ii) $\reg(S/x S) \leq \reg(S) = r-1$. Again  by 
Claim 2 this yields
$$\reg(K^1(S/xS)) <  4(\reg(S/xS))^{n-1} - 4(\reg(S/xS))^{n-2} \leq  4 r^{n-1} - 4 r^{n-2} = \alpha  .$$
Hence we can apply Claim 4 with $i=1$ and $\alpha > 0$ to get
$$|P_{\overline{K^2}}(-1)| < \frac{1}{2} r^{n-3} (r + 4 r^{n-1} - 4 r^{n-2}
+ 3)^4.$$
By Claim 3  and  (\ref{EC3}) we obtain
$$\begin{array}{ll}
\Delta_1 &< \Delta_0 + r^{n-1} +\frac{1}{2} r^{n-3} (r + 4 r^{n-1} - 4 r^{n-2}
+ 3)^4\\
& < 3r^{n-1} - 2 r^{n-2} + \frac{1}{2} r^{n-3} ( 4 r^{n-1} - 4r)^4 \ \ \text{(since}\ n\geq 
4,\ r\geq 2)\\
& < 3r^{n-1} - 2 r^{n-2} + \frac{1}{2} r^{n+1} 4^4 (  r^{4(n-2)} - 4) \\
& < \frac{1}{2}(2 r)^{n(n+1)} - r^n -n \ \text{(since}\ n\geq 4).
\end{array}$$
Thus the first inequality is proven. 

Furthermore, by the inequalities at the beginning of the proof, we can use Claim 1 to get
$$\reg(K^2/xK^2) \leq \alpha +2 = 4 r^{n-1} - 4 r^{n-2} + 2 .$$ 
Hence, by Lemma \ref{C2} and Lemma
\ref{A7} (ii), this  implies
$$\begin{array}{ll}
\reg(K^2) & \leq \max\{ \reg(K^2/xK^2) ,\  [2(1+ \Delta_1)]^2 -2 \}\\
 & \leq \max\{ 4 r^{n-1} - 4 r^{n-2} + 2 ,\ [2(1+ \frac{1}{2}(2 r)^{n(n+1)} - r^n -n )]^2 
-2\}\\
& \leq (2r)^{2n(n+1)} - 2r^n - 2n. 
\end{array}$$
This is the second inequality of the claim. \hfill $\square$
\end{pf}

\vskip0.3cm

\noindent {\it CLAIM 6}. Assume that $1 \leq i < d-1$. Then
$$\Delta_i <\frac{1}{2}(2 \reg(I))^{n\cdots (n+i)2^{\frac{i(i-1)}{2}}} - (\reg I)^n -n,$$
and
$$\reg(K^{i+1}) <(2 \reg(I))^{n\cdots (n+i)2^{\frac{i(i+1)}{2}}}- 2(\reg I)^n -2n.$$

\begin{pf} We do induction on $i$. The case $i=1$ is Claim 5. Let $i\geq 2$ and let $x$ 
be a generic linear element.  
By the  induction hypothesis we have
$$\reg K^i< \beta - 2r^n -2n,$$
where
$$\beta := (2 r)^{n\cdots (n+i-1)2^{\frac{i(i-1)}{2}}}.$$
Since $\reg(S/xS) \leq \reg(S) = r-1$, the induction hypothesis also gives
$$ \begin{array}{ll}
\reg K^i(S/xS) &< [2(\reg S/xS) + 1)]^{n\cdots (n+i-1)2^{\frac{i(i-1)}{2}}} - 2(\reg S/xS +1)^n - 2n\\
& \leq (2 r)^{n\cdots (n+i-1)2^{\frac{i(i-1)}{2}}} - 2r^n - 2n =
 \beta - 2r^n -2n.
\end{array}$$
Applying Claim 3, Claim 4 (with $\alpha := \beta - 2r^n -2n$) and the induction
hypothesis on $\Delta_{i-1}$, we get
$$\begin{array}{ll}
\Delta_i & < i\Delta_{i-1} + (\reg I)^{n-1} + |P_{\overline{K^{i+1}}}(-i)| \\
& \leq i(\frac{1}{2}\beta - r^n -n) + r^{n-1} + \frac{1}{2}
r^{n-i-2}(r + \beta
- 2r^n -2n + 2i + 1)^{2i+2}\\
 & < \frac{1}{2} \beta^{2i+3} - r^n -n\\
 &< \frac{1}{2} \beta^{n+i} - r^n -n  \ \text{(since}\ n\geq i+3) \\
 & = \frac{1}{2}(2 r)^{n\cdots (n+i)2^{\frac{i(i-1)}{2}}} - r^n -n.
\end{array}$$
Furthermore, by recalling the inequalities at the beginning of the induction step, 
we can use Claim 1 to have
$$\reg(K^{i+1}/xK^{i+1}) \leq \beta - 2r^n -2n +2 .$$ 
By Lemma \ref{C2} and Lemma
\ref{A7} (ii),  this now  implies that 
$$\begin{array}{ll}
\reg(K^{i+1}) & \leq \max\{ \reg(K^{i+1}/xK^{i+1}) ,\  [2(1+ \Delta_i)]^{2^i} -2  \}\\
 & \leq \max\{ \beta - 2r^n -2n + 2 ,\ [2(1+ \frac{1}{2} \beta^{n+i} - r^n -n  )]^{2^i}   
-2\}\\
& <  \beta^{(n+i)2^i}   - 2r^n - 2n \\
& = (2 r)^{n\cdots (n+i)2^{\frac{i(i+1)}{2}}} - 2r^n - 2n . 
\end{array}$$
Claim 6 is thus completely proven. \hfill $\square$
\end{pf}

The cases $i=1$ and $2 \leq i \leq d-1$ of Theorem \ref{C1} were proved in Claim 2 and
 Claim 6, respectively. To finish the proof of Theorem \ref{C1},  we only have  to show 
the following stronger bound

\vskip0.3cm

\noindent {\it CLAIM 7}. Let $d\geq 2$. Then
$$\reg(K^d) < (2 \reg(I))^{n\cdots (n+d-2)2^{\frac{(d-1)(d-2)}{2}}}- 2(\reg I)^n -2n 
+2.$$

\begin{pf}  Let
$$\beta(d) = (2 r)^{n\cdots (n+d-2)2^{\frac{(d-1)(d-2)}{2}}} - 2r^n -2n.$$ 
We will prove by  induction on $d$ that $ \reg(K^d) < \beta(d)+2$. 

Let $d=2$.  This case was considered in \cite{HaH}, Theorem  2.9 and the bound there 
is much  smaller. For the convenience of the reader we give here a direct proof of the 
weaker bound: $ \reg(K^2) <  \beta(2)+2$. Let  $x$ be,  as usual,  a generic 
linear element. Since $\reg (S/xS) \leq \reg( S) = r-1$ by Claim 2, both $\reg(K^1)$ and 
$\reg(K^1 (S/xS))$  are less than $4r^{n-1} - 4 r^{n-2}$. Since $n\geq 3$, 
$$ 4r^{n-1} - 4 r^{n-2} < \beta(2) = (2r)^n - 2r^n - 2n.$$
By (\ref{EB5b})  we then get  $ \reg(K^2) <  \beta(2)+2$. 

Let $d\geq 3$. By Claim 6,  
$$\reg(K^{d-1}) < \beta(d).$$
 Since $\dim S/xS = d-1$, by the  induction hypothesis the following  holds
 $$\begin{array}{ll}
\reg(K^{d-1}(S/xS))  & <  [ 2(\reg(S/xS) +1)]^{n\cdots (n+d-3)2^{\frac{(d-2)(d-3)}{2}}}\\
& \hskip1cm  - 2(\reg(S/xS) +1)^n - 2n + 2\\
   & <   (2 r)^{n\cdots (n+d-2)2^{\frac{(d-1)(d-2)}{2}}} - 2r^n -2n =  \beta (d).
\end{array}$$
 Hence,  again by (\ref{EB5b}),  we  get 
$\reg(K^d) < \beta(d) +2 ,$
 as required. \hfill $\square$
\end{pf}

\begin{rem}\label{C4} {\rm Assume that $S$ is a generalized Cohen-Macaulay ring, i.e.
all modules $K^i(S),\ i<d$, are of finite length. In this case $q^i_S(t) =0$
for all $i\leq d-2$, and the proof of Theorem \ref{C1} will be substantially
simplified. It gives
$$\reg(K^i(S)) < [2^{i+1}( (\reg I)^{n-1} - (\reg
I)^{n-2})]^{2^{i-1}}.$$ It is still a huge number. We don't know whether one
can give a linear bound even in this case.}
\end{rem}

The bounds in Theorem \ref{C1} are huge. However, this theorem demonstrates
that the Castelnuovo-Mumford regularity $\reg(S)$ also controls the behavior
of local cohomology modules in negative components. To understand better this
phenomenon, let us state some consequences. The first corollary is formulated
in the spirit of \cite{BMM}, Theorem 4.8.

\begin{cor}\label{C5} Denote by $\Hcal_{n,i,r}$ the set of  numerical
functions $h:\ \Zset \rightarrow \Zset$ such that there exists a homogeneous
ideal $I \subset R= k[x_1,...,x_n]$ satisfying the following conditions
\begin{itemize}
\item[(i)] $\reg I \leq r$,
\item[(ii)] $\ell(H^i_{\mm}(R/I)_t) = h(t)$ for all $t\in \Zset$.
\end{itemize}
Then for fixed numbers $n,i,r$ the set $\Hcal_{n,i,r}$ has only finitely many
elements. 
\end{cor}

\begin{pf} Note that $h(t) = 0$ for all $t \geq r$. By Theorem
\ref{C1}, $\reg(K^i(S))$ is bounded by a number $f(n,r)$ depending on $n$ and
$r$. By Lemma \ref{C3}, for each $t$ with $ - (f(n,r) + n) \leq t \leq r$, the
value $h(t) = \ell(H^i_{\mm}(R/I)_t)$ is also bounded by a function $g(n,r)$.
This implies that there are only finitely many choices of the initial values
of $h(t)$. Since $P_{K^i}(t)= \ell(H^i_{\mm}(R/I)_{-t})$ in $n$ points $t=
f(n,r)+1,...,f(n,r) + n$, and the degree of $P_{K^i}(t)$ is less than $n$, the
number of possible polynomials $P_{K^i}(t)$ is finite. Moreover $h(t) =
 P_{K^i}(- t)$ for all $t< -f(n,r)$. These statements  together imply the
 finiteness of the set $\Hcal_{n,i,r}$. \hfill $\square$
\end{pf}

Assume that $k$ is an algebraically closed field. A famous result of Kleiman
states that there exists only a finite number of Hilbert functions associated
to reduced and equi-dimensional $k$-algebras $S$ such that $\deg(S) \leq e$
and $\dim(S) = d$. In a recent paper \cite{Hoa} the first author was able to
extend this result to all reduced algebras. Recall that
$$\adeg S =  \sum_{\pfrak \in \text{Ass}(S)} \ell(H^0_{\mm_{\pfrak}}(S_{\pfrak}))
 e(S/\pfrak)$$
is called arithmetic degree of $S$ (see \cite{BM}, Definition 3.4 or
\cite{Va2}, Definition 9.1.3). The arithmetic degree agrees with $\deg(S)$ if
and only if $S$ is equi-dimensional. Inspired by Kleiman's result we formulate
the following two corollaries.

\begin{cor}\label{C6} Denote by $\Hcal_{d,i,a}^*$ the set of all numerical
functions $h:\ \Zset \rightarrow \Zset$ for which there exists a reduced
$k$-algebra $S$ such that $\adeg(S) \leq a$, $\dim (S) = d$ and
$\ell(H^i_{\mm}(S)_t) = h(t)$ for all $t\in \Zset$. Assume that $k$ is an
algebraically closed field. Then for fixed numbers $d,i,a$ the set
$\Hcal_{d,i,a}^*$ is finite.
\end{cor}

\begin{pf} Under the assumption, by \cite{Hoa}, Theorem 1.5, $\reg(S)$ is
bounded by $f(d,a)$, and by \cite{Hoa}, Lemma 5.2, $n$ is bounded by $g(d,a)$
too. Hence the assertion follows from Corollary \ref{C5}. \hfill $\square$
\end{pf}

\begin{cor}\label{C7} Denote by $\Hcal'_{n,i,\delta}$ the set of all numerical
functions $h:\ \Zset \rightarrow \Zset$ such that there exists a homogeneous
ideal $I \subset R= k[x_1,...,x_n]$ satisfying the following conditions
\begin{itemize}
\item[(i)] $I$ is generated by forms of degrees at most $\delta$,
\item[(ii)] $\ell(H^i_{\mm}(R/I)_t) = h(t)$ for all $t\in \Zset$.
\end{itemize}
Then for fixed numbers $n,i,\delta$ the set $\Hcal'_{n,i,\delta}$ has only
finitely many elements.
\end{cor}

\begin{pf} Under the assumption, by \cite{Hoa}, Theorem 2.1 (see also
\cite{BM}, Proposition 3.8), $\reg(S)$ is bounded by $f(n,\delta)$. Hence the
assertion follows from Corollary \ref{C5}. \hfill $\square$
\end{pf}

\section{Examples}\smallskip \label{D}

We believe that there should be a much better bound for $\reg(K^i(S))$ in
terms of $\reg(S)$ than the one given in the previous section. In this section
we show this for some particular cases.

\vskip0.3cm

 1. If $S=R/I$ is the coordinate ring of a smooth projective variety over a
 field of characteristic zero, then M. Chardin and B. Ulrich (\cite{CU},
 Theorem 1.3) showed that
 $$\reg(K^d(S)) = d.$$
 This is a consequence of Kodaira's vanishing theorem.

\vskip0.3cm

2. Assume that $M$ is a generalized Cohen-Macaulay module, i.e. all modules
$K^i(M),\ i<d$, are of finite length. In the general case there is no known
good bound for $\reg(K^i(M))$ (see Remark \ref{C4}). However, there is a good
one in terms of the annihilators of $K^i(M)$ (see \cite{HM}, Proposition 2.4
and Corollary 2.5). We recall here a nice case.

A module $M$ is called Buchsbaum module if the difference $\ell(M/\qfrak M) -
e(\qfrak, M)$ between the length and the multiplicity is a constant, when
$\qfrak$ runs over all homogeneous parameter ideals of $M$. In this case $\mm
K^i(M) = 0$ for all $i< d$.  Proposition 2.4 (i) in \cite{HM} states that if
$M$ is a Buchsbaum module, then
\begin{equation}\label{ED1}
\reg(K^i(M))  \leq i-\bg(M),\ i\leq d.
\end{equation}

3. Assume that $I$ is a monomial ideal and that $S=R/I$ is a generalized
Cohen-Macaulay ring. Then by \cite{T}, Proposition 1 we have
\begin{equation}\label{ED2}
\reg(K^i(S)) = \ed (K^i(S)) \leq 0\ \text{for}\ i<d.
\end{equation}
 We also get
\begin{prop}\label{D1} Assume that $I$ is a monomial ideal and that $S=R/I$ is a generalized
Cohen-Macaulay ring. Then
$$\reg(K^d(S)) \leq d.$$
\end{prop}

\begin{pf} Since $S$ is a generalized Cohen-Macaulay ring, we have  by \cite{Sc1}, Corollary 3.1.3
 the  following isomorphisms
$$H^{d+1-i}_{\mm}(K^d) \cong K^i,$$
 for all $2\leq i < d$, and
there is an exact sequence of graded modules
$$ 0 \rightarrow K^1 \rightarrow H^{d}_{\mm}(K^d) \rightarrow \Hom (S,k)
\rightarrow K^0 \rightarrow 0.$$ We also have $\dept(K^d)\geq \min\{2,
\dim(S)\}$. Combining this  with (\ref{ED2}) implies the assertion.
\end{pf}

4. In some cases when $M$ is not necessarily a generalized Cohen-Macaulay
module,  good bounds can still be found for $\reg(K^d(M))$ (see \cite{HaH}).
In order to extend these results to all $\reg(K^i(M))$, let us recall some
definitions.

For an integer $0 \leq i \leq d$, let $M^i$ denote the largest graded
submodule of $M$ such that $\dim M^i \leq i$. Let $M^{-1} = 0$. The increasing
filtration
$$0 = M^{-1} \subseteq M^0 \subseteq \cdots \subseteq M^d = M$$
is called the dimension filtration of $M$. This filtration is well-defined and
unique. We put
$$\Mcal^i = M^i/M^{i-1} \ \ \text{for\ all }\ 0\leq i\leq d.$$
Note that $\Mcal^i$ is either zero or of dimension $i$. A module $M$ is called
a sequentially Cohen-Macaulay (sequentially Buchsbaum) module if each module
$\Mcal^i$ is either zero or a Cohen-Macaulay (Buchsbaum, respectively). The
notion of a sequentially Cohen-Macaulay module was introduced by R. Stanley
(see, e.g., \cite{HS}).

\begin{prop}\label{D2} (i) If $M$ is a sequentially Cohen-Macaulay
module, then for all $i\leq d$ we have
$$\reg(K^i(M))\leq i- \bg(M).$$
(ii)  If $M$ is a sequentially Buchsbaum module, then for all $i\leq d$ we
have
$$\reg(K^i(M))\leq i +1 - \bg(M).$$
\end{prop}
\begin{pf} (i) Under the assumption, $K^i(M) \cong K^i(\Mcal^i)$ by \cite{Sc2},
Lemma 5.2. If $\Mcal^i=0$, then  there is nothing to prove. Otherwise,
$\Mcal^i$ is a Cohen-Macaulay module of dimension $i$. By (\ref{EB2b}) we have
$$\reg(K^i(M))= \reg(K^i(\Mcal^i)) = i - \bg(\Mcal^i) \leq i - \bg(M).$$

(ii) We do induction on $d$. If $d=1$ then $M$ is a sequentially Cohen-Macaulay module. 
Hence the assertion holds true by (i).  Let $d\geq 2$.
By Lemma \ref{B7bn} we may assume that $i>0$.  
The case $i=d$ is \cite{HaH}, Proposition 2.2. Let $1\leq i<d$.  The exact
sequence
$$0 \rightarrow M^{d-1} \rightarrow M \rightarrow \Mcal^d \rightarrow 0$$
gives the long exact sequence of cohomology
$$ K^i(\Mcal^d) \overset{\varphi}\rightarrow K^i(M) \overset{\psi}\rightarrow
K^i(M^{d-1}) \overset{\chi}\rightarrow K^{i-1}(\Mcal^d).$$ This breaks up into
two short exact sequences
$$\begin{array}{l}
0 \rightarrow \text{Im}\, \varphi \rightarrow K^i(M) \rightarrow
\text{Im}\,\psi \rightarrow 0, \\
0 \rightarrow \text{Im}\, \psi \rightarrow K^i(M^{d-1}) \rightarrow
\text{Im}\,\chi \rightarrow 0.
\end{array}$$
Note that $\dim(\Mcal^d) = d$, and that by the assumption $\Mcal^d$ is a
Buchsbaum module. Hence $K^i(\Mcal^d)$ and $K^{i-1}(\Mcal^d)$ are modules of
finite length. Since $\bg(\Mcal^d) \geq \bg(M)$,  we have by (\ref{ED1})
$$\reg(\text{Im}\, \varphi) \leq \reg(K^i(\Mcal^d)) \leq i- \bg(\Mcal^d) \leq
i - \bg(M),$$
and
$$\reg(\text{Im}\, \chi) \leq \reg(K^{i-1}(\Mcal^d)) \leq i- 1 - \bg(M).$$
Using Lemma \ref{A3} and the above two short exact sequences we obtain
$$\begin{array}{ll}
\reg(K^i(M)) & \leq \max\{\reg(\text{Im}\, \varphi) , \ \reg(\text{Im}\, \psi)
\} \\
& \leq \max\{i - \bg(M) ,\ \reg(K^i(M^{d-1})), \ \reg(\text{Im}\, \chi) + 1
\} \\
& = \max\{i - \bg(M) ,\ \reg(K^i(M^{d-1}))\}.
\end{array}$$
Since $M^{d-1}$ is also a sequentially Buchsbaum module (of dimension at most
$d-1$),  we know by the induction hypothesis that
$$ \reg(K^i(M^{d-1})) \leq i - \bg(M^{d-1}) \leq  i-\bg(M).$$
Consequently, $\reg(K^i(M))\leq  i-\bg(M).$ \hfill $\square$
\end{pf}

Let $\gin(I)$ denote the generic ideal of $I$ with respect to a term order. It
is a so-called Borel-fixed ideal, and by \cite{HS}, Theorem 2.2, $R/\gin(I)$
is a sequentially Cohen-Macaulay ring. Hence we get:

\begin{cor}\label{D3} For an arbitrary homogeneous ideal $I\subset R$ we
have
$$\reg(K^i(R/ \gin I)) \leq i.$$
\end{cor}

Unfortunately we cannot use this result to bound $\reg(K^i(R/I))$. The reason
is the following. We always have
 $$\ell(K^i(R/\gin I)_j) \geq \ell(K^i(R/I)_j) \ \ \text{for\ all}\ j\in \Zset.$$
 It is well-known that many invariants  increase by passing from $I$ to $\gin(I)$,
 but remain unchanged if one takes the generic initial ideal $\Gin(I)$ with respect 
to the reverse
lexicographic order. So, if   the equality
$$\ell(K^i(R/\Gin I)_j) = \ell(K^i(R/I)_j)$$
would hold for all $j\in \Zset$, then using Corollary \ref{D3} one would get a good 
bound for
$ri(K^i(R/I))$. From that, by the method of Section \ref{C}, one would get a
good bound for $\reg(K^i(R/I))$. Unfortunately, this is almost impossible.
Namely, J. Herzog and E. Sbarra (\cite{HS}, Theorem 3.1) proved that
$$\ell(K^i(R/\Gin I)_j) = \ell(K^i(R/I)_j)   ,$$
for all $i \leq d$ and all $j\in \Zset$ if and only if $R/I$ is itself a
sequentially Cohen-Macaulay ring. But the latter case was settled in Proposition \ref{D2} (i).

\vskip0.5cm \noindent {\bf Acknowledgement}. The first named author would like
to thank the Academy of Finland for financial support and the University of
Helsinki for hospitality during the preparation of this article.

\end{document}